\documentclass[twosided,reqno]{amsart}

\usepackage{amsmath}
\usepackage{amssymb}
\usepackage{amsthm}

\theoremstyle{theorem}
\newtheorem{theorem}{\scshape Theorem }[section]

\newtheorem{corollary}[theorem]{\scshape Corollary}

\theoremstyle{definition}

\newcommand{\ma}{\mathbb}
\newcommand{\be}{\begin{equation}}
\newcommand{\ee}{\end{equation}}
\newcommand{\ben}{\begin{equation*}}
\newcommand{\een}{\end{equation*}}
\newcommand{\fa}{\frac}
\newcommand{\la}{\label}
\newcommand{\Z}{\mathbb{Z}_{p}}
\newcommand{\m}{d\mu_{-1}}
\newcommand{\M}{d\mu_{-1}(x_1)\ldots d\mu_{-1}(x_k)}
\newcommand{\E}{E_{n}}
\newcommand{\ek}{E_{n}^{(k)}}
\newcommand{\so}{S_{1}}
\newcommand{\st}{S_{2}}
\newcommand{\ty}{\infty}
\newcommand{\I}{\int_{\mathbb{Z}_{p}}\cdots \int_{\mathbb{Z}_{p}}}
\newcommand{\U}{\sum_{n=0}^{\infty}}
\newcommand{\bt}{\begin{theorem}}
\newcommand{\et}{\end{theorem}}
\newcommand{\bi}{\binom}
\newcommand{\al}{\alpha}
\newcommand{\lp}{\left(}
\newcommand{\rp}{\right)}

\newcommand{\CK}{\widehat{Ch}_n^{(k)}}
\newcommand{\ck}{Ch_n^{(k)}}

\numberwithin{equation}{section}

\begin{document}

\title{Higher-order  Changhee numbers and polynomials  }

\author{Dae San Kim}
\address{Department of Mathematics, Sogang University, Seoul 121-742, Republic of Korea.}
\email{dskim@sogang.ac.kr}

\author{Taekyun Kim}
\address{Department of Mathematics, Kwangwoon University, Seoul 139-701, Republic of Korea}
\email{tkkim@kw.ac.kr}

\maketitle
\begin{abstract}
In this paper, we consider the higher-order Changhee numbers and polynomials which are derived from the fermionic $p$-adic integral on $\Z$ and give some relations between higher-order Changhee polynomials and special polynomials.
\end{abstract}

\section{Introduction}

As is well known, the Euler polynomials of order $\al(\in \ma{N})$ are defined by the generating function to be

\be\la{1}
\lp \fa{2}{e^{t}+1}\rp^{\alpha}  e^{xt}=\U \E^{(\al)}(x)\fa{t^{n}}{n!},~(\textrm{see [1-16]}).
\ee
When $x=0$, $\E^{(\al)}=\E^{(\al)}(0)$ are called the Euler numbers of order $\alpha$.\\
The Stirling number of the first kind is defined by

\be\la{2}
(x)_{n}=\sum_{l=0}^{n}\so(n,l)x^{l},~(n \in \ma{Z}_{\geq 0}),~ (\textrm{see [5,6,7]}),
\ee
where $(x)_{n}=x(x-1)\cdots (x-n+1)$.\\
The Stirling number of the second kind is also defined by the generating function to be

\be\la{3}
(e^t-1)^n=n!\sum_{l=n}^{\ty}\st(l,n)\fa{t^{l}}{l!},~(n \in \ma{Z}_{\geq 0}).
\ee

Let $p$ be an odd prime number.  Throughout this paper, $\Z$, $\ma{Q}_{p}$ and $\ma{C}_{p}$ will denote the ring of $p$-adic integers, the field of $p$-adic numbers and the completion of algebraic closure of $\ma{Q}_{p}$.  The $p$-adic norm $|\cdot|_{p}$ is normalized as $|p|_{p}=\fa{1}{p}$.  Let $C(\Z)$ be the space of continuous functions on $\Z$. For $f \in C(\Z)$, the fermionic $p$-adic integral on $\Z$ is defined by Kim to be
\be\la{4}
I_{-1}(f)=\int_{\mathbb{Z}_{p}} f(x)\m(x)=\lim_{N \to \ty}\sum_{x=0}^{p^{N}-1}f(x)(-1)^{x},~(\textrm{see \cite{09}}).
\ee

For $f_{1}(x)=f(x+1)$, we have

\be\la{5}
I_{-1}(f_{1})+I_{-1}(f)=2f(0).
\ee

As is well-known, the Changhee polynomials are defined by the generating function to be

\ben
\fa{2}{t+2}(1+t)^{x}=\sum_{n=0}^{\ty}Ch_n(x)\fa{t^{n}}{n!},~(\textrm{see [6,8]}).
\een
When $x=0$, $Ch_n=Ch_n(0)$ are called the Changhee numbers. In this paper, we consider the higher-order Changhee numbers and polynomials which are derived from the multivariate fermionic $p$-adic integral on $\Z$ and give some relations between higher-order Changhee polynomials and special polynomials.

\section{Higher-order Changhee  polynomials}

For $k\in \ma{N}$, let us define the Changhee numbers of the first kind with order $k$ as follows:

\be\la{6}
\ck=\I(x_1+\cdots +x_k)_n\M,
\ee
where $n$ is a nonnegative integer.\\
From (\ref{6}), we can derive the generating function of $\ck$ as follows:
\be\la{7}
\begin{split}
\U\ck \fa{t^{n}}{n!}&=\I\U\bi{x_1+\cdots +x_k}{n}t^n\M\\
&=\I(1+t)^{x_1+\cdots +x_k}\M.
\end{split}
\ee

By (\ref{5}), we easily see that

\be\la{8}
\I(1+t)^{x_1+\cdots +x_k}\M=\lp\fa{2}{2+t}\rp^k.
\ee

From (\ref{7}) and (\ref{8}), we have

\be\la{9}
\U\ck \fa{t^{n}}{n!}=\lp\fa{2}{2+t}\rp^k.
\ee

It is easy to show that

\be\la{10}
\lp\fa{2}{2+t}\rp^k=\U\lp\sum_{l_1+\cdots+l_k=n}\bi{n}{l_1,\cdots,l_k}Ch_{l_1}\cdots Ch_{l_k} \rp \fa{t^{n}}{n!}.
\ee

Thus, by (\ref{9}) and (\ref{10}), we get
\be\la{11}
\ck=\sum_{l_1+\cdots+l_k=n}\bi{n}{l_1,\cdots,l_k}Ch_{l_1}\cdots Ch_{l_k}.
\ee

It is not difficult to show that

\be\la{12}
\lp\fa{2}{2+t}\rp^k=\U\lp -\fa{1}{2}\rp^n n!\bi{k+n-1}{n}\fa{t^{n}}{n!}.
\ee

From (\ref{9}) and (\ref{12}), we have
\be\la{13}
\begin{split}
2^n\ck&=(-1)^n n!\bi{n+k-1}{n}=(-1)^n(k+n-1)_n\\
&=(-1)^n\sum_{l=0}^{n}\so(n,l)(k+n-1)^l.
\end{split}
\ee

Therefore, by (\ref{13}), we obtain the following theorem.

\bt\la{t1}
For $n \geq 0$, we have
\ben
\ck=\lp -\fa{1}{2}\rp^n \sum_{l=0}^{n}\so(n,l)(k+n-1)^l.
\een
\et

By (\ref{6}), we get

\be\la{14}
\begin{split}
\ck&=\I (x_1+\cdots + x_k)_n\M\\
&=\sum_{l=0}^{n}\so(n,l)\I(x_1+\cdots +x_k)^l\M.
\end{split}
\ee
Now, we observe that

\be\la{15}
\I e^{(x_1+\cdots x_k)t}\M=\lp\fa{2}{e^t+1}\rp^k=\U \ek \fa{t^{n}}{n!}.
\ee

By (\ref{14}) and (\ref{15}), we get

\be\la{16}
\ck=\sum_{l=0}^{n}\so(n,l)E_l^{(k)}.
\ee

Therefore, by (\ref{16}), we obtain the following theorem.

\bt\la{t2}
For $n\geq 0$, we have
\ben\la{16}
\ck=\sum_{l=0}^{n}\so(n,l)E_l^{(k)}.
\een
\et

Replacing $t$ by $e^t-1$ in (\ref{9}), we get

\be\la{17}
\U \ck \fa{(e^t-1)^{n}}{n!}=\lp\fa{2}{e^t+1}\rp^k=\sum_{m=0}^{\ty}E_m^{(k)}\fa{t^{m}}{m!},
\ee
and

\be\la{18}
\U \ck \fa{(e^t-1)^{n}}{n!}=\sum_{m=0}^{\ty}\lp \sum_{n=0}^{m}\ck \st(m,n)\rp \fa{t^{m}}{m!}.
\ee

Therefore, by (\ref{17}) and (\ref{18}), we obtain the following theorem.

\bt\la{t3}
For $n \geq 0$, we have
\ben
E_m^{(k)}=\sum_{n=0}^{m}\ck \st(m,n).
\een
\et

Now, we consider the higher-order Changhee polynomials of the first kind as follows:
\be\la{19}
\ck(x)=\I(x_1+\cdots x_k+x)_n\M.
\ee

By (\ref{19}), we get

\be\la{20}
\U \ck(x)\fa{t^{n}}{n!}=\I(1+t)^{x_1+\cdots x_k+x}\M=\lp\fa{2}{2+t}\rp^k (1+t)^x.
\ee

From (\ref{9}), we have
\be\la{21}
\lp\fa{2}{2+t}\rp^k(1+t)^x=\U\lp \sum_{m=0}^{n}\bi{n}{m}(x)_m Ch_{n-m}^{(k)}\rp \fa{t^{n}}{n!}.
\ee

By (\ref{20}) and (\ref{21}), we get

\be\la{22}
\ck(x)=\sum_{m=0}^{n}\bi{x}{m}\fa{n!}{(n-m)!} Ch_{n-m}^{(k)}=\sum_{m=0}^{n}\bi{x}{n-m}\fa{n!}{m!} Ch_{m}^{(k)}.
\ee

From (\ref{19}), we have
\be\la{23}
\begin{split}
\ck(x)&=\I(x_1+\cdots x_k+x)_n\M\\
&=\sum_{l=0}^{n}\so(n,l)E_l^{(k)}(x).
\end{split}
\ee

Therefore, by (\ref{23}), we obtain the following corollary.

\begin{corollary}\la{c4}
For $n\geq 0$, we have
\ben
\ck(x)=\sum_{l=0}^{n}\so(n,l)E_l^{(k)}(x).
\een
\end{corollary}

In (\ref{20}), by replacing $t$ by $e^t-1$, we get

\be\la{24}
\U\ck(x)\fa{(e^t-1)^{n}}{n!}=\lp\fa{2}{e^t+1}\rp^ke^{tx}=\sum_{m=0}^{\ty}E_m^{(k)}(x)\fa{t^{m}}{m!},
\ee
and

\be\la{25}
\U\ck(x)\fa{1}{n!}(e^t-1)^{n}=\sum_{m=0}^{\ty}\lp \sum_{n=0}^{m}\ck(x)\st(m,n)\rp \fa{t^{m}}{m!}.
\ee

Therefore, by  (\ref{24}) and  (\ref{25}), we obtain the following theorem.

\bt\la{t5}
For $m\geq 0$, we have
\ben
E_m^{(k)}(x)=\sum_{n=0}^{m}\ck(x)\st(m,n).
\een
\et

The rising factorial is defined by

\be\la{26}
(x)^{(n)}=x(x+1)\cdots (x+n-1)=(-1)^n(-x)_n.
\ee
Here, we define the Changhee numbers of the second kind with order $k(\in \ma{N})$ as follows:
\be\la{27}
\CK=\I(-x_1- \cdots -x_k)_n\M.
\ee

Thus, by (\ref{27}), we get
\be\la{28}
\begin{split}
\CK&=\sum_{l=0}^{n}(-1)^l\so(n,l)\I(x_1+ \cdots x_k)^l \M\\
&=\sum_{l=0}^{n}(-1)^l\so(n,l)E_l^{(k)}.
\end{split}
\ee

The generating function of $\CK$ is given by
\be\la{29}
\begin{split}
\U\CK \fa{t^{n}}{n!}&=\I(1+t)^{-x_1- \cdots -x_k}\M\\
&=\lp\fa{2}{2+t}\rp^k(1+t)^k.
\end{split}
\ee

Now, we observe that
\be\la{30}
\lp\fa{2}{2+t}\rp^k(1+t)^k=\U\lp \sum_{m=0}^{n}\bi{k}{m}Ch_{n-m}^{(k)}\fa{n!}{(n-m)!}\rp\fa{t^n}{n!}.
\ee

Thus, by (\ref{29}) and (\ref{30}), we get

\be\la{31}
\CK = \sum_{m=0}^{n} m!\bi{k}{m}\bi{n}{m}Ch_{n-m}^{(k)}.
\ee

Therefore, by (\ref{31}), we obtain the following theorem.

\bt\la{t6}
For $n \geq 0$, we have
\ben
\CK = \sum_{m=0}^{n} m!\bi{k}{m}\bi{n}{m}Ch_{n-m}^{(k)}.
\een
\et

In  (\ref{29}), by replacing $t$ by $e^t-1$, we get

\be\la{32}
\U\CK \fa{(e^t-1)^{n}}{n!}=\lp\fa{2}{e^t+1}\rp^k e^{tk}=\sum_{m=0}^{\ty}E_m^{(k)}(k)\fa{t^{m}}{m!},
\ee
and
\be\la{33}
\U\CK \fa{(e^t-1)^{n}}{n!}=\sum_{m=0}^{\ty}\lp \sum_{n=0}^{m}\CK\st(m,n)\rp \fa{t^{m}}{m!}.
\ee

Therefore, by  (\ref{32}) and  (\ref{33}), we obtain the following theorem.

\bt\la{t7}
For $m \geq 0$, we have
\ben
E_m^{(k)}(k)=\sum_{n=0}^{m}\CK\st(m,n).
\een
\et

Now, we consider the Changhee polynomials of the second kind with order $k(\in \ma{N})$ as follows:

\be\la{34}
\CK(x)=\I(-x_1- \cdots -x_k+x)_n\M.
\ee

From (\ref{29}) and (\ref{34}), we have
\be\la{35}
\begin{split}
\U \CK(x)\fa{t^{n}}{n!}&=\I(1+t)^{-x_1- \cdots -x_k+x}\M\\
&=(1+t)^{x+k}\lp\fa{2}{2+t}\rp^{k}.
\end{split}
\ee

We observe that

\be\la{36}
\lp\fa{2}{2+t}\rp^{k}(1+t)^{x+k}=\U\lp \sum_{m=0}^{n} m!\bi{x}{m}\bi{n}{m}Ch_{n-m}^{(k)}\rp \fa{t^{n}}{n!}.
\ee

Thus, by (\ref{35}) and (\ref{36}), we obtain the following theorem.

\bt\la{t8}
For $m \geq 0$, we have
\ben
\CK(x)=\sum_{m=0}^{n} m!\bi{x}{m}\bi{n}{m}Ch_{n-m}^{(k)}.
\een
\et

From (\ref{34}), we have
\be\la{37}
\begin{split}
\CK(x)&=\sum_{l=0}^{n}\so(n,l)(-1)^l\I (x_1+ \cdots +x_k-x)^l \M\\
&=\sum_{l=0}^{n}\so(n,l)(-1)^l E_l^{(k)}(-x).
\end{split}
\ee

In (\ref{35}), by replacing $t$ by $e^t-1$, we get

\be\la{38}
\U\CK(x)\fa{(e^t-1)^{n}}{n!}=e^{(x+k)t}\lp\fa{2}{e^t+1}\rp^k=\sum_{m=0}^{\ty}E_m^{(k)}(x+k)\fa{t^{m}}{m!},
\ee
and
\be\la{39}
\U\CK(x)\fa{1}{n!}(e^t-1)^{n}=\sum_{m=0}^{\ty}\lp \sum_{n=0}^{m}\CK\st(m,n)\rp \fa{t^{m}}{m!}.
\ee

Therefore, by (\ref{38}) and (\ref{39}), we obtain the following theorem.

\bt\la{t9}
For $m \geq 0$, we have
\ben
E_m^{(k)}(x+k)=\sum_{n=0}^{m}\CK\st(m,n).
\een
\et

Now, we observe that
\be\la{40}
\begin{split}
(-1)^n\fa{\CK(x)}{n!}&=(-1)^n\I\bi{-(x_1+ \cdots +x_k)+x}{n}\M\\
&=\I\bi{x_1+ \cdots +x_k-x+n-1}{n}\M\\
&=\sum_{m=0}^{n}\bi{n-1}{n-m}\I\bi{x_1+ \cdots +x_k-x}{m}\M\\
&=\sum_{m=0}^{n}\fa{\bi{n-1}{n-m}}{m!}Ch_m^{(k)}(-x)=\sum_{m=1}^{n}\fa{\bi{n-1}{n-m}}{m!}Ch_m^{(k)}(-x).
\end{split}
\ee
Therefore, by (\ref{40}), we obtain the following theorem.

\bt\la{t10}
For $n\in \ma{N}$, we have
\ben
(-1)^n\fa{\CK(x)}{n!}=\sum_{m=1}^{n}\fa{\bi{n-1}{n-m}}{m!}Ch_m^{(k)}(-x).
\een
\et

By (\ref{19}), we get

\be\la{41}
\begin{split}
(-1)^n\fa{\ck(x)}{n!}&=(-1)^n\I\bi{x_1+ \cdots +x_k+x}{n}\M\\
&=\I\bi{-x_1-x_2- \cdots -x_k-x+n-1}{n}\M\\
&=\sum_{m=0}^{n}\fa{\bi{n-1}{n-m}}{m!}\widehat{Ch}_m^{(k)}(-x)=\sum_{m=1}^{n}\fa{\bi{n-1}{n-m}}{m!}\widehat{Ch}_m^{(k)}(-x).
\end{split}
\ee

Therefore, by (\ref{41}), we obtain the following theorem.

\bt\la{t11}
For $n\in \ma{N}$, we have
\ben
(-1)^n\fa{\CK(x)}{n!}=\sum_{m=1}^{n}\fa{\bi{n-1}{n-m}}{m!}\widehat{Ch}_m^{(k)}(-x).
\een
\et


\end{document}